\documentclass[12pt]{article}
\usepackage[utf8]{inputenc}
\usepackage{latexsym}     
\usepackage{amsfonts}
\usepackage{amsthm}
\usepackage{amsmath}
\usepackage{amssymb}
\usepackage{graphicx}

\theoremstyle{remark}

\def\tm{{\pitchfork}}
\def\mt{{\emptyset}}

\def\bfA{{\bf A}}

\def\d{\partial}

\def\Z{\mathbb{Z}}

\title{The combinatorial transverse intersection algebra}
\author{Daniel An, Ruth Lawrence, Dennis Sullivan}
\date{03/24}

\begin{document}

\maketitle
\begin{abstract}

  {\scriptsize
  
  This paper constructs (with challenging obstacles) on the three torus with its cubical decomposition: 
   
   Firstly, a  combinatorial  graded intersection algebra (graded by the codimension)  which is commutative and  associative  defined by transversality   on the usual chains  which are in  general position. This, (with extra elements added)  on the
   entire $h$-cubulated three torus whose differential satisfies the product rule and which agrees  with the set theoretic  intersection product appropriately weighted.  The construction  is  characterized  given these properties (see Comprehensive Theorem below). The challenge  is to minimally adjoin  infinitesimal elements when the geometric elements have glancing but transversal intersections weighted  in such a way that the associativity (and commutativity)  is not destroyed  and  the Leibniz product rule for the boundary operator is restored. 
   
   Secondly, there is  a $2h$ subcomplex introduced in [3] and discussed further in [2] which shares the above good properties when ideal elements are added  AND  which  also has a star bijection between degree zero and degree three and between one and degree two. This  is introduced for the purposes of computations of 3D fluid motion, incompressible, with or without viscosity [1]. For the latter  purposes one only needs the good properties in dimensions zero, one and two, where the situation is a bit better.

   It is a new feature that the three good properties are respected by the crumbling   chain mappings from coarse to finer subdivisions. The star operator  does not cooperate  with crumbling and is the sole reason in this discrete approximation for the Kolmogorov cascade to finer scales.
   
   }
\end{abstract}
\section{Aim and background}
Consider a cubical lattice with lattice spacing $h$ in three dimensions. The usual associated chain complex has four non-trivial chain spaces, in dimensions 0,1,2 and 3 which have basis elements (see Figure 1) which are points, elemental edges of length $h$ parallel to one of the three axes, elemental plaquettes which are squares of edge length $h$ parallel to one of the three coordinate planes, and elemental cubes of edge length $h$, respectively. 
\[
\includegraphics[width=.6\textwidth]{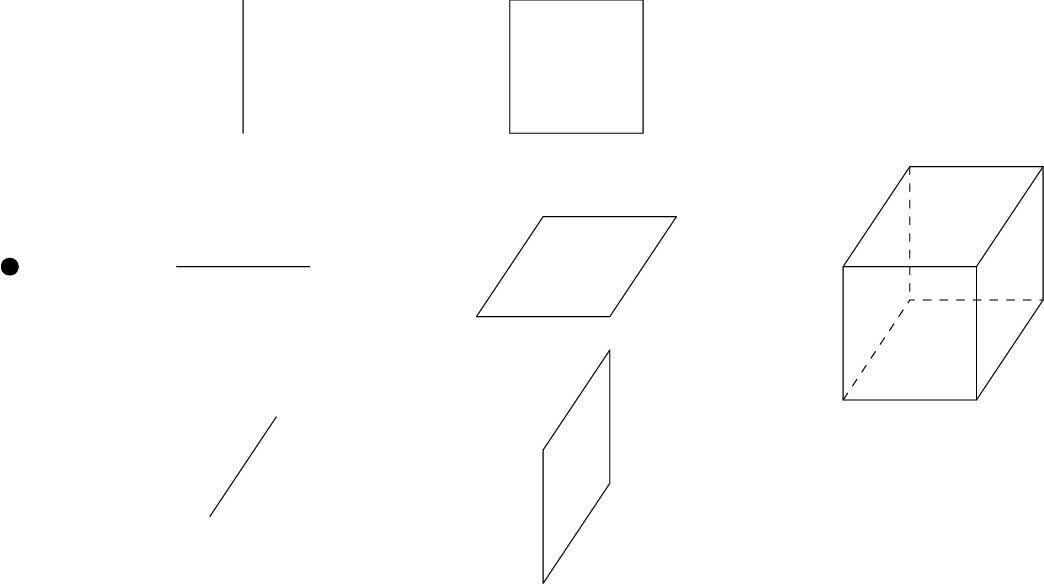}
\]\[\hbox{\small\sl Figure 1: The $h$ complex}
\]

\medskip
\noindent The boundary map is defined by the geometric boundary (with orientations).
 
 Inside this chain complex we will also consider for applications the $2h$ subcomplex which has basis consisting of all possible similar cells but with edge length $2h$, as in Figure 2.

\[ 
\includegraphics[width=.5\textwidth]{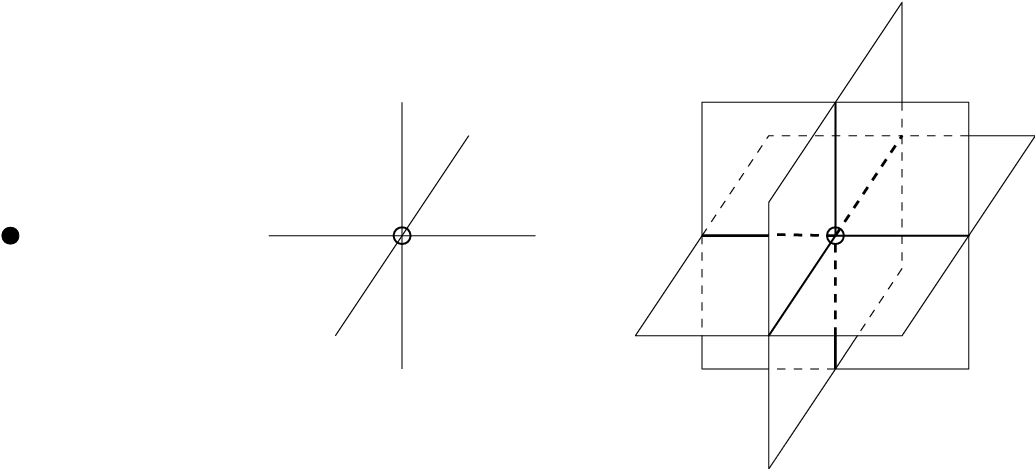}\hskip2em
\includegraphics[width=.22\textwidth]{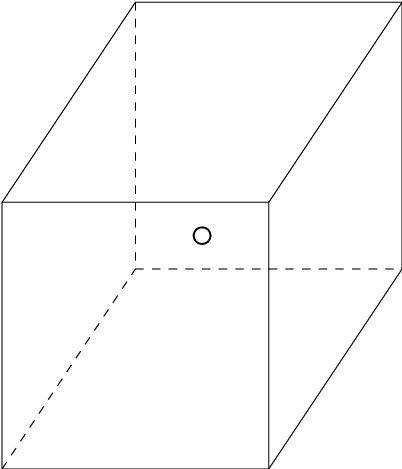}
\]
\[\hbox{\small\sl Figure 2: The $2h$ complex}\]

\medskip \noindent In order to consider the $2h$ chain spaces as subspaces of the $h$ chain spaces, we consider the larger cells as subdivided into cells with edge length $h$, so that $2h$ edges are sums of pairs of connected parallel $h$-edges, $2h$-squares are sums of four connected coplanar $h$-squares fanning a square of edge length $2h$, and $2h$-cubes are sums of eight $h$-cubes.

The $2h$ subcomplex has nice properties:  each vertex is the barycenter of (1,3,3,1) $2h$ cells of dimensions 0,1,2,3 as in   Figure 2  and these form an  exterior algebra structure at each vertex. The entire $2h$  chain complex with boundary operator of degree $-1$  has a graded commutative associative (with the sign determined by the codimensions) intersection algebra structure  with   degrees $i,j$ going  to  degree $i+j -3$. Plus there is a natural star duality relating degrees one and two and  relating degrees zero and three. 

These two chain complexes, the h cubical decomposition and the subcomplex of overlapping two h cells  is a first approximation of the  dual intersection  geometric picture of the  exterior algebra structure on  differential forms with the geometric  boundary operator of degree $-1$  being in  Hom duality with  the  picture of the exterior  derivative on forms \cite{S}. 

 The caveat   is that the $\d$ operator of degree -1   does not satisfy  the product rule for this first approximation of the geometric product.  This  discrepancy   has been  treated  firstly, by an infinity algebraic structure   \cite{LRS} but then the problem arose that such PDEs as Euler or Navier-Stokes were not  obvious to write  in that enlarged context. 
 
 \bigskip 
 The point of this paper is to  further develop this combinatorial intersection  product  to  substantially  restore the product rule for the boundary operator with  a  better approximation,  truer to  the geometric  intersection product.   This successful  modification will be seen to also apply to the original $h$-complex where Poincar\'e duality is also largely restored in terms of the intersection pairing.
 The product will respect  the  intersections that are geometrically transverse and in  general position.

A parallelopiped whose edges are parallel to the coordinate axes and whose vertices lie in the lattice will be called a {\it cuboid}. A cuboid can be subdivided into  $h$-cubes, see Figure 3(a). Similarly in lower dimensions, rectangles and sticks on the lattice, of arbitrary edge lengths (multiples of $h$) can be subdivided into $h$-cells; such objects we will call {\it cuboidal cells} on the lattice. Thus a cuboidal cell in the lattice is geometrically a Cartesian product of singletons and intervals (all of whose delimiters lie in $h\Z$). Replacing one or more of the intervals defining a cuboidal cell by singletons at one of the endpoints of the corresponding intervals will lead to geometric objects which are cuboidal cells of lower dimension, namely faces, edges and vertices of the original cuboid; we will call them {\it generalised faces} of the cuboidal cell, see Figure 3(b). Each cuboid cell is a sum of basis (cubical) cells of the $h$-complex, and in the case of cuboids all of whose edge lengths are even, it can be expressed as a sum of cells of the $2h$-complex. 
\[\hskip4em\includegraphics[width=.3\textwidth]{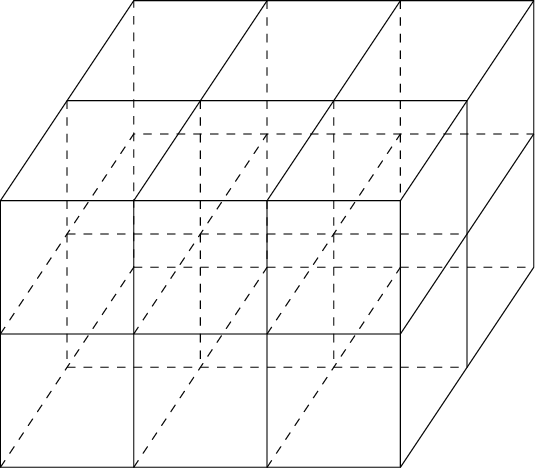}\hskip6em
\includegraphics[width=.3\textwidth]{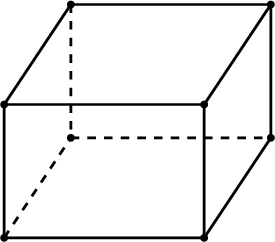}\]
\[\hbox{\small\sl Figure 3:(a) Subdivision of a cuboid into cubes (b) Generalised faces of a cuboid}\]

\medskip Now consider a pair of cuboidal cells. The geometric intersection of two cuboidal cells is considered to be {\it transverse} if the set-wise intersection of the closed cells is non-empty while their tangent spaces generate the entire (three-dimensional) tangent space; for example, two intersecting lines are not transverse in three-dimensions. We say that two cuboidal cells are in {\it general position} if they intersect transversely, and in addition whenever we replace either or both cuboid by one of its generalised faces, all such pairs of cuboidal cells are either disjoint or have transverse intersection. The possible configurations of pairs of  cuboid $2h$-cells in three dimensions in general position are shown in Figure 4. 
\[\includegraphics[width=.46\textwidth]{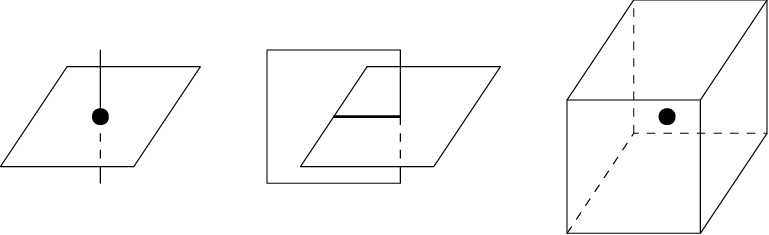}\hskip2em
\includegraphics[width=.46\textwidth]{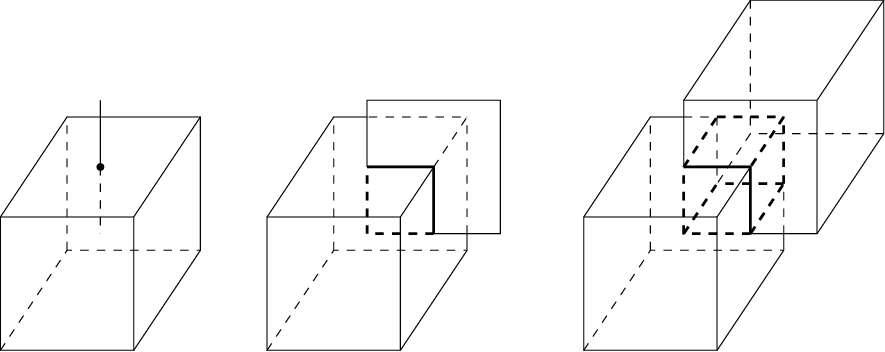}
\]
\[\hbox{\small\sl Figure 4: Intersections in the $2h$ complex in general position}\]

\medskip\noindent Observe that the geometric intersection of cuboidal cells in general position is also a cuboidal cell whose dimension is the sum of the dimensions of the initial cells minus three; equivalently the codimension of the intersection is the sum.

The possible types of configurations of $h$-cells in three dimensions which intersect transversally but are {\sl not} in general position are shown in Figure 5.
\[\includegraphics[width=.19\textwidth]{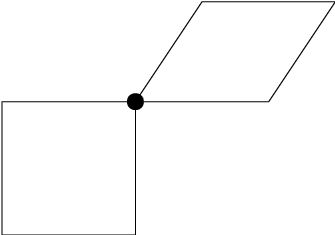}\hskip3ex
\includegraphics[width=.11\textwidth]{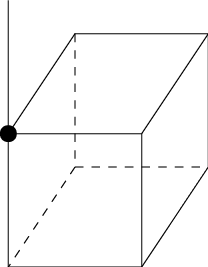}\hskip4ex
\includegraphics[width=.4\textwidth]{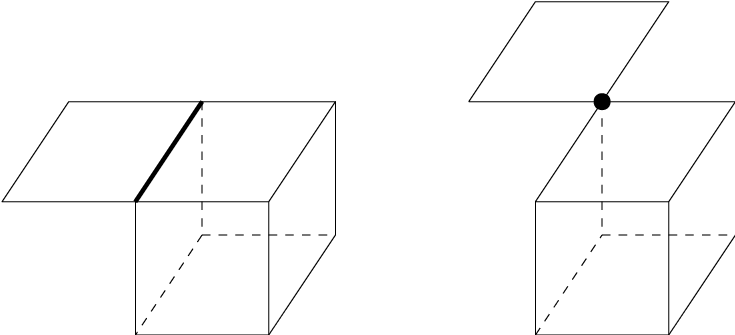}\]
\[\includegraphics[width=.78\textwidth]{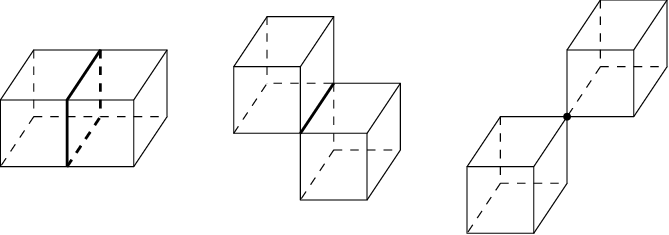}\]
\[\hbox{\small\sl Figure 5: Cuboid intersections suggest the ideal elements added to the $h$-complex}\]

\medskip\noindent In the first example, the intersection `should' be of dimension $2+2-3=1$ but is geometrically a point. There are five more types of ideal elements to be added as suggested by the remaining examples (see Figure 9).

Dealing with this, which will result in restoring the product rule for $\d$ plus other benefits,   is the subject of our discussion.

 We now  describe a finite algebraic structure extending the usual chain complex of a cubical lattice whose multiplication is described by geometric intersection in the case of intersecting cuboids in general position, and which is commutative, associative and satisfies (as closely as possible) the product rule with respect to the boundary operator.
The resulting chain complex we will call the {\sl enlarged $h$-cubical chain complex with infinitesimals}.  

In order to obtain a finite-dimensional algebra, it is necessary to use a periodic lattice. In order that the $2h$ cell subcomplex will compute the correct rational homology groups of the three-torus, it is necessary and sufficient that the period of the lattice be odd in each of the three directions.  Otherwise there is  a copy of the homology for mod two  equivalence class  of vertices ie 1,2,4,8 copies if three, two, one or none of the periods are odd( exercise for the  interested reader).

%\vfill
\section{Comprehensive Theorem}
 
    The intersection geometry  announced  above can be completely  worked out algebraically in dimension one. Then one can tensor copies of these any number of times to obtain the following theorem in dimension three and in other dimensions. The proofs for dimension one and three  follow  directly and further  dimensions and  requisite  supplements  are in the   Addendum.

\medskip\medskip\noindent
{\bf Theorem}  By adding  finitely many infinitesimal elements suggested by figure 5 , there exists a chain complex in three dimensions which is an enlargement $EC_*$ of the usual triply periodic (even or odd independently in all directions) cubical complex $C_*$, on which there is a locally defined\footnote{By {\sl locally defined} we mean that the set theoertic support of the intersection of a pair of basis elements is the geometric intersection of the supports of the basis elements.} `transverse intersection' product $EC_*\otimes{}EC_*\to{}EC_*$ satisfying
\begin{itemize}
    \item[(A)] graded commutativity (the grading of cells is by codimension\footnote{This may be surprising, but seems natural when it is understood that this intersection product is a geometric picture of the exterior product of differential forms, where a plaquette corresponds to a 1-form}), that is, $a\cdot{}b=(-1)^{c_ac_b}b\cdot{}a$
    \item[(B)] associativity, $(a\cdot{}b)c=a(b\cdot{}c)$
    \item[(C)] the product rule for the boundary operator holds on the original cubical complex $C_*$, that is $\d(a\cdot{}b)=(\d{a})\cdot{}b+(-1)^{c_a}a\cdot(\d{b})$,  and the  product on the cubical complex  contains infinitesimal elements. That is, the original subcomplex is not closed under the intersection product.
    \item[(D)] the product is invariant under the symmetries of the cubical lattice
    \item[(E)] the product is non-zero on pairs of cuboidal cells precisely when they are  transverse, that is their closures have non-trivial set theoretic intersection  plus their tangent spaces generate all of three-dimensional space
    \item[(F)] the product agrees  on cuboid cells with the natural intersection product (signed geometric intersection of closed cells) in the case of  cuboidal cells  which are in general position
    \item[(G)] the augmentation map in degree zero which sums all the coefficients applied to any    product to degree zero  defines a pairing $\langle{}a,b\rangle$ such that (1) (Frobenius associativity) $\langle{}a\cdot{}b,c\rangle=\langle{}a,b\cdot{}c\rangle$ holds; (2) when the periods in all directions are odd, the pairing is non-degenerate on the original $C_*$ subcomplex  (which recall  has no infinitesimal elements);
    \item[(H)] there is a subalgebra $FC_*$ of $EC_*$ containing the cubical complex $C_*$ in dimensions 0,1,2 on which the product rule holds without reservation\footnote{Note: 
 We use the ``$2h$'' subcomplex of $FC_*$, which has a star operator, to study fluids \cite{AnKwon}.}.
     \item[(I)] There is a unique minimal extension of $C^*$ with a product satisfying (A)-(F).
    \item[(J)] The natural (crumbling) chain mapping from a coarse periodic cubical complex to a finer periodic cubical complex commutes with the intersection product satisying (A)-(F).
\end{itemize} 

{\bf Remark}
Actually the product on the three-dimensional $h$ cubical complex in dimensions up to two, with ideal sticks $FC_*$ adjoined, is unique satisfying commutativity, associativity and the product rule assuming locality and that it extends the natural geometric intersection on all pairs of cuboids in general position.  For the detailed statement See \S5.

\section{One-dimensional case}
To understand the previous construction and to extend it further to all dimensions it is useful to treat the one-dimensional case first.
\medskip

\noindent{\bf Theorem} There exists a locally defined (transverse intersection) product on the one-dimensional chain complex $C_*$ with odd period enlarged to $EC_*$ by infinitesimal elements, which satisfies
\begin{itemize}
    \item[(A)] graded commutativity (the grading is the codimension)
    \item[(B)] associativity
    \item[(C)] the boundary satisfies the product rule on the subcomplex $C_*$
    \item[(D)] the product is respected by the symmetries of the lattice
    \item[(E)] the product is non-zero on pairs of cells precisely when their closures have non-trivial geometric intersection and their tangent spaces generate (one-dimensional) space
    \item[(F)] the product agrees with the natural intersection product (geometric intersection of closed cells) in the case of pairs of cells (points, intervals of arbitrary length in the lattice) in general position (that is, having distinct points and endpoints)
    \item[(G)] the augmentation of the product defines a pairing $\langle{}a,b\rangle$ such that (1) (Frobenius associativity) $\langle{}a\cdot{}b,c\rangle=\langle{}a,b\cdot{}c\rangle$ holds; (2) the pairing is non-degenerate on the subcomplex without infinitesimal elements (odd period)
\end{itemize}
There is a unique product with these properties.
Additionally the natural (crumbling) chain mapping from a coarse  odd period one-dimensional complex to a finer odd period one-dimensional complex commutes with the intersection product.
\begin{proof}
To demonstrate uniqueness, we use the given properties to identify the form of product in terms of certain parameters which are then forced to satisfy certain equations. Solving the equations reveals a single solution (up to scaling of the infinitesimal elements) which is then verified to indeed satisfy all the properties and thus demonstrates also existence.

\[\includegraphics[width=.4\textwidth]{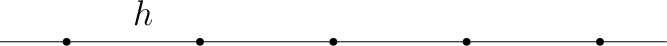}\hskip4em
\includegraphics[width=.2\textwidth]{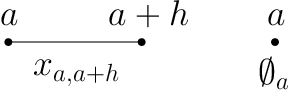}
\]
\[\hbox{\small\sl Figure 6: (a) One-dimensional $h$-lattice (b) Generating cells}\]

\medskip
In this one-dimensional case, the initial $h$-complex has dimension one cells generated by sticks of length $h$ (say from $a\in{}h\Z$ to $a+h$) denoted $x_{a,a+h}$, with dimension zero cells (points) denoted by $\mt_a$, $a\in{}h\Z$.  The boundary map is given by
\[
\d(x_{ab})=\mt_b-\mt_a
\]
More general cuboidal cells are simply longer sticks,
$x_{ab}$ for $a,b\in{}h\Z$, $a<b$, which may be decomposed as
\[x_{ab}=x_{a,a+h}+x_{a+h,a+2h}+\cdots+x_{b-h,b}\eqno{(1)}\]
The product preserves the codimension grading so that the non-trivial products are point$\tm$stick (which gives a codimension $1+0$ object, that is a point) and stick$\tm$stick (which gives a codimension $0+0$ object, that is a stick). Since all non-trivial intersections involve at least one even (zero) graded object, graded commutativity (A) is just plain commutativity in the one-dimensional case.

\[\includegraphics[width=.7\textwidth]{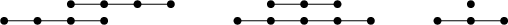}\]
\[\hbox{\small\sl Figure 7: Intersections of cuboids in general position in one dimension}\]

\medskip
Transverse intersection occurs between points and sticks and between pairs of sticks, so long as there is a non-empty geometric intersection. This generates infinitesimal objects from the intersection of sticks whose geometric intersection is a point as in the middle diagram of Figure 8; denote these by $x_{a,a}$ for $a\in{}h\Z$.
\[\includegraphics[width=.7\textwidth]{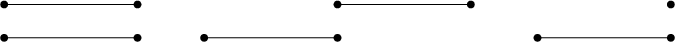}\]
\[\hbox{\small\sl Figure 8: Transverse intersections of $h$-cells not in general position in 1-d}\]

\medskip\noindent
The relevant intersections of basis elements in the enlarged complex are thus
\begin{align*}
    \hbox{point$\tm$stick}:&\quad\mt_a\tm{}x_{a,a+h},\quad \mt_a\tm{}x_{a-h,a},\quad \mt_a\tm{}x_{a,a}\\
    \hbox{stick$\tm$stick}:&\quad{}x_{a-h,a}\tm{}x_{a,a+h},\quad x_{a,a+h}\tm{}x_{a,a+h}
\end{align*}
By (A), (E), these are all non-zero and are the only non-trivial cases (up to reversal of order). Locality limits the possible terms which can be involved in the product, to be those of the correct dimension and supported in the geometric intersection; in the case of a point intersection this leaves only a non-zero scalar parameter. Since the symmetries of the one-dimensional lattice are translation by $h$ and reflection $x\longmapsto{}-x$, thus locality and (D) ensure that for some non-zero constants $s,t$ (independent of $a$)
\[\mt_a\tm{}x_{a,a+h}=s\mt_a,\qquad \mt_a\tm{}x_{a-h,a}=s\mt_a,\qquad \mt_a\tm{}x_{a,a}=t\mt_a\]
For the intersection of sticks, the result is one-dimensional and so locality and (D) ensure that for some scalars $\alpha,\beta,\gamma,\delta,\epsilon$ independent of $a$,
\begin{align*}x_{a-h,a}\tm{}x_{a,a+h}&=\alpha{}x_{a,a},\qquad 
x_{a,a+h}\tm{}x_{a,a+h}=\beta{}x_{a,a}+\gamma{}x_{a,a+h}+\beta{}x_{a+h,a+h},\\
x_{a,a}\tm{}x_{a,a+h}&=\delta{}x_{a.a},\qquad{}x_{a,a}\tm{}x_{a,a}=\epsilon{}x_{a,a}
\end{align*}
where $\alpha,\delta,\epsilon\not=0$ while $(\beta,\gamma)\not=(0,0)$. Furthermore, we can choose the scaling of the infinitesimal element $x_{a,a}$ so that $\alpha=1$, meaning that we choose $x_{a,a}\equiv{}x_{a-h,a}\tm{}x_{a,a+h}$.

By (C), $\d(x\tm{}y)=(\d{}x)\tm{}y+(-1)^{|x|}x\tm\d{y}$ when $x,y$ are in the initial $h$-complex (that is, without infinitesimal elements). This is only non-trivial when $x$ and $y$ are both sticks. By linearity, it is enough to verify this for basis elements, the only two non-trivial cases being
\begin{align*}
    \d(x_{a-h,a}\tm{}x_{a,a+h}):&\ \d(\alpha\mt_a)=(\mt_a-\mt_{a-h})\tm{}x_{a,a+h}+x_{a-h,a}\tm(\mt_{a+h}-\mt_a)\\
    \d(x_{a,a+h}\tm{}x_{a,a+h}):&\ \d(\beta\mt_a+\gamma{}x_{a,a+h}+\beta\mt_{a+h}))=2(\mt_{a+h}-\mt_a)\tm{}x_{a,a+h}&
\end{align*}
The first equality reduces to $0=s\mt_a-s\mt_a$ while the second gives $\gamma(\mt_{a+h}-\mt_a)=2s(\mt_{a+h}-\mt_a)$, that is $\gamma=2s$.

Intersection of cuboids in general position occurs so long as in addition there are no coincidences of (end)points, namely the possible intersections in general position of cuboids in one-dimension are
\begin{align*}
   x_{a,b}\tm{}\mt_c&\qquad\hbox{for $a<c<b$}\\
   x_{a,b}\tm{}x_{c,d}&\qquad\hbox{for $a,b,c,d$ distinct and $[a,b]\cap[c,d]\not=\mt$}
\end{align*}
By (F), we know that in this case the intersection must be the geometric intersection,
\begin{align*}
   x_{a,b}\tm{}\mt_c&=\mt_c\qquad\,\,\,\,\hbox{for $a<c<b$}\\
   x_{a,b}\tm{}x_{c,d}&=x_{c,b}\qquad\hbox{for $a<c<b<d$}\\
   x_{a,b}\tm{}x_{c,d}&=x_{c,d}\qquad\hbox{for $a<c<d<b$}
\end{align*}
Using the linearity of the product $\tm$ and the decomposition (1) of long sticks as a sum of sticks of length $h$, these relations give
\begin{align*}
   2s\mt_c&=\mt_c\\
   (\alpha+\beta)(x_{c,c}+\cdots+x_{b,b})+\gamma{}x_{c,b}&=x_{c,b}\qquad\hbox{for $a<c<b<d$}\\
   (\alpha+\beta)(x_{c,c}+\cdots+x_{d,d})+\gamma{}x_{c,d}&=x_{c,d}\qquad\hbox{for $a<c<d<b$}
\end{align*}
where the intermediate terms are $2x_{p,p}$ for all intermediate lattice points. These equalities require $s=\frac12$, $\gamma=1$ and $\alpha+\beta=0$. Note that this ensures $\gamma=2s$ as required for (C). So far we have found that these conditions are guaranteed by, and guarantee (A),(C),(D),(E) and (F). This leaves only (B) and (G).

To ensure associativity (B), it is enough to verify it on triples of basis elements in the $h$-complex (plus infinitesimals). Since multiplication takes a pair of dimensions $i,j$ to an object of dimension $i+j-1$, the only relevant dimensions for triples are $1,1,0$ and $1,1,1$. The only triples $x,y,z$ of basis elements (up to symmetries) in which at least one of the three possible products $x(yz)$, $(xy)z$, $(xz)y$ is non-zero are
\begin{align*}
    1,1,0:&\qquad\includegraphics[width=.4\textwidth]{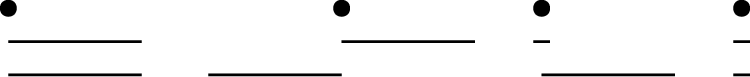}\\
    &\\
    1,1,1:&\qquad\includegraphics[width=.7\textwidth]{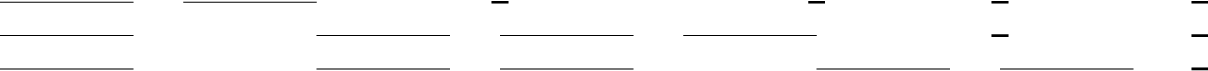}
\end{align*}
Associativity is tautologous for a triple of identical elements. The rest of the triples give conditions for associativity as follows.
\begin{align*}
\includegraphics[width=.1\textwidth]{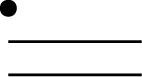}:\qquad{}\beta{}t+\gamma{}s=s^2\\[1ex]
\includegraphics[width=.18\textwidth]{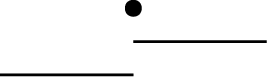}:\qquad{}\alpha{}t=s^2\\[1ex]
\includegraphics[width=.09\textwidth]{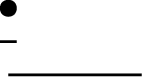}:\qquad{}\delta{}t=st\\[1ex]
\includegraphics[width=.01\textwidth]{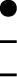}:\qquad{}\epsilon{}t=t^2
\end{align*}

\begin{align*}
\includegraphics[width=.18\textwidth]{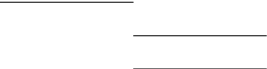}:\qquad{}\beta\delta+\gamma\alpha=\alpha\delta\\[1ex]
\includegraphics[width=.1\textwidth]{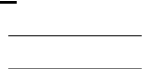}:\qquad{}\beta\epsilon+\gamma\delta=\delta^2\\[1ex]
\includegraphics[width=.18\textwidth]{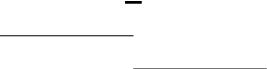}:\qquad{}\alpha\epsilon=\delta^2\\[1ex]
\includegraphics[width=.1\textwidth]{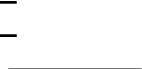}:\qquad{}\delta\epsilon=\epsilon\delta
\end{align*}
These equations have the solution
\[
s=\delta=\frac12,\quad \gamma=1,\quad \epsilon=t,\quad \alpha=\frac1{4t},\quad \beta=-\frac1{4t}
\]
Having chosen above the scaling of the infinitesimal element so that $\alpha=1$, this requires $t=\frac14$ and so the unique multiplication satisfying (A)-(F) is defined on basis elements by 
\begin{align*}
\mt_a\tm{}x_{a,a+h}&=\frac12\mt_a,\qquad \mt_a\tm{}x_{a-h,a}=\frac12\mt_a,\qquad \mt_a\tm{}x_{a,a}=\frac14\mt_a\\
x_{a-h,a}\tm{}x_{a,a+h}&=x_{a,a},\qquad 
x_{a,a+h}\tm{}x_{a,a+h}=-x_{a,a}+x_{a,a+h}-x_{a+h,a+h},\\
x_{a,a}\tm{}x_{a,a+h}&=\frac12x_{a.a},\qquad{}x_{a,a}\tm{}x_{a-h,a}=\frac12x_{a.a},\qquad{}x_{a,a}\tm{}x_{a,a}=\frac14x_{a,a}
\end{align*}
The coefficients here may be understood naturally in terms of the probabilistic interpretation \cite{AnKwon} whereby each generator is considered with endpoints which are not exactly on the lattice, but rather shifted slightly according to a probability distribution with support $\epsilon<<h$, so that all intersections are in general position with probability 1.

Finally we deal with (G). The augmentation is the `counting points' map $\#:C_0\to\Z$ which maps $\mt_a\longmapsto1$. The augmentation of the product defines a symmetric pairing by 
\[\langle{}a,b\rangle=\#(a\tm{}b)\]  
This can be extended to a map $\langle\>,\>\rangle:EC_*\times{}EC_*\to\Z$ (zero except on the parts $C_0\times{}EC_1$ and $EC_1\times{}C_0$) with the property
\[\langle{}a\tm{}b,c\rangle=\#((a\tm{}b)\tm{}c)=\#(a\tm(b\tm{}c))=\langle{}a,b\tm{}c\rangle\]
Meanwhile the pairing restricted to $C_0\times{}C_1$ with respect to the basis $\{\mt_a\}$ for $C_0$ and $\{x_{a,a+h}\}$ for $C_1$ has matrix elements
\[A_{a,b}=\langle\mt_{ah},x_{bh,(b+1)h}\rangle=\frac12\delta_{a,b}+\frac12\delta_{a,b+1}\]
which is a cyclic $N\times{}N$ matrix with entries $\frac12,0,\ldots,0,\frac12$, where $N$ is the period of the lattice. The eigenvalues of this matrix are $\frac12+\frac12\omega^j$ where $\omega=\exp{\frac{2\pi{}i}{N}}$ is a primitive $N-th$ root of unity. Thus the determinant of the matrix of the pairing is 
\[\det{\bfA}=\prod\limits_{j=0}^{N-1}\left(\frac12+\frac12\omega^j\right)\]
This is non-zero (and hence the pairing is non-degenerate) so long as $\omega^j\not=-1$ for all $j$, that is so long as $N$ is odd. This completes the proof of (G).

Finally, consider the crumbling map $\iota:EC_*(h)\longrightarrow{}EC_*(h')$ where $h'=\frac{h}{k}$ for $k$ odd. This transforms a one-dimensional lattice of period $N$ (odd) to a one-dimensional lattice of period $Nk$ (still odd), will preserve points and infinitesimal sticks while breaking elemental $h$ sticks in the coarse lattice into a sum of $k$ elemental sticks in the fine lattice,
\begin{align*}
    \mt_a&\longmapsto\mt_a\\
    x_{a,a}&\longmapsto{}x_{a,a}\\
    x_{a,a+h}&\longmapsto{}x_{a,a+h'}+x_{a+h',a+2h'}+\cdots+x_{a+(k-1)h',a+kh'}
\end{align*}
It commutes with the multiplication just found, the only instance which isn't immediate being the case of $x_{a,a+h}\tm{}x_{a,a+h}$ which evaluates as a sum of $k^2$ terms, of which $3k-2$ are non-vanishing (diagonal and first off-diagonals) leading to
\begin{align*}
&\sum\limits_{i=1}^kx_{a+(i-1)h',a+ih'}\tm{}x_{a+(i-1)h',a+ih'}
+2\sum\limits_{i=1}^{k-1}x_{a+(i-1)h',a+ih'}\tm{}x_{a+ih',a+(i+1)h'}\\
&=\sum\limits_{i=1}^k\Big(x_{a+(i-1)h',a+ih'}-x_{a+(i-1)h',a+(i-1)h'}-x_{a+ih',a+ih'}\Big)
+2\sum\limits_{i=1}^{k-1}x_{a+ih',a+ih'}\\
&=\sum\limits_{i=1}^k\Big(x_{a+(i-1)h',a+ih'}\Big)-x_{a,a}-x_{a+kh',a+kh'}
\end{align*}
as needed.\end{proof}

\section{Proof of comprehensive theorem}
The cubic tensor power of the one-dimensional algebra gives an example of a transverse multiplication algebra for the three-dimensional lattice satisfying (A)-(G) and commuting with the crumbling map. The elements of this complex $EC_*$ will be those of the original cubical complex $C_*$ along with additional ideal elements coming from tensor products of one or more infinitesimal elements from the one-dimensional complex. Indeed there are six types of ideal elements arising from these tensor products, see Figure 9.

\break\[\includegraphics[width=.7\textwidth]{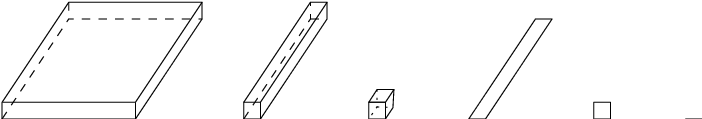}\]
\[\hbox{\small\sl Figure 9: Ideal elements in $EC_*$}\]
For (H), note that the only transverse intersections of cuboids of dimension $\leq2$ which are basis elements in the $h$-complex $C_*$ and which are not in general position (thus generating infinitesimal elements) are of the type
\[\includegraphics[width=.2\textwidth]{nonstrictintersect22}\]
Let $FC_*$ be the subcomplex of $EC_*$ generated by $C_{\leq2}$ along with only infinitesimal (one-dimensional) sticks. Such generators will be labelled by their direction and location $x_{pp}$, $y_{pp}$, $z_{pp}$ where $p$ is a point in the lattice. This subcomplex is closed under multiplication, forming a subalgebra of $EC_*$. It remains only to verify that the product rule holds on the whole of $FC_*$. From the one-dimensional theorem, we know that the product rule 
\[\d(a\tm{}b)=(\d{a})\tm{}b+(-1)^{c_a}a\tm(\d{b})\]
holds for $a,b\in{}C_*$. If both of $a,b$ are infinitesimal sticks then the product rule holds trivially (all terms vanish from dimensional considerations). If $a$ is an infinitesimal stick and $b\in{}C_{\leq2}$ then again the product rule holds trivially (all three terms vanish from dimensional considerations).

\section{Uniqueness of 3-d commutative associative Leibniz product on  the $F$ subalgebra}
In the previous section we constructed an intersection product on the three-dimensional enlarged complex $EC_*$ satisfying (A)-(G) by taking the tensor product of three copies of the one-dimensional intersection product. In \S3, we showed that the one-dimensional intersection product satisfying (A)-(G) is unique; this does not imply that the three-dimensional intersection product is unique. Also in the last section, motivated by the  precise needs of computing the 3D Navier Stokes Fluid  evolution , we showed the existence in three-dimensions of a subalgebra $FC_*$ of $EC_*$ on which the product is commutative, associative {\sl and} Leibniz, spanned by $C_{\leq2}$ and infinitesimal (one-dimensional) sticks. In this section we show that the product on $FC_*$ is unique satisfying these properties, explicitly  one has the following theorem.  

\medskip\noindent{\bf Theorem} On the subcomplex $FC_*$ of $EC_*$ spanned by the dimension $\leq2$ part of the periodic cubical complex complex $C_*$ and infinitesimal (one-dimensional) sticks, there exists a unique (up to rescaling of infinitesimal elements) locally defined `transverse intersection' product $FC_*\otimes{}FC_*\to{}FC_*$  satisfying
\begin{itemize}
    \item[(A)] graded commutativity, $a\cdot{}b=(-1)^{c_ac_b}b\cdot{}a$
    \item[(B)] associativity, $(a\cdot{}b)c=a(b\cdot{}c)$
    \item[(C)] the product rule,  $\d(a\cdot{}b)=(\d{a})\cdot{}b+(-1)^{c_a}a\cdot(\d{b})$
    \item[(D)] the product is invariant under the symmetries of the cubical lattice
    \item[(E)] the product is non-zero on pairs of cuboidal cells precisely when they are  transverse, that is their closures have non-trivial geometric intersection and their tangent spaces generate all of three-dimensional space
    \item[(F)] the product agrees with the natural intersection product (signed geometric intersection of closed cells) in the case of cuboidal cells in general position
    \item[(G)] the augmentation of the product defines a pairing $\langle{}a,b\rangle$ such that (1) (Frobenius associativity) $\langle{}a\cdot{}b,c\rangle=\langle{}a,b\cdot{}c\rangle$ holds; (2) when the periods are odd in all directions, the pairing is non-degenerate on the original  cubical $C_*$ subcomplex (which has no  infinitesimal elements).
\end{itemize}
Additionally the natural (crumbling) chain mapping from a coarse periodic cubical complex to a finer periodic cubical complex commutes with the intersection product.
\begin{proof}
Existence was proved in \S{4}, so it remains to prove uniqueness. The only non-trivial products in $FC_*$ are $C_2\times{}C_2\to{}FC_1$ and $C_2\times{}FC_1\to{}C_0$.
Since both have codimension one, a graded symmetric transverse intersection of squares is anti-symmetric. It also changes sign under change of orientation of either square. Intersections of a pair of rectangles in general position look like Figure 10(a); we will use the case of two $2h$-squares in this configuration. Property (F) requires that their product be the natural geometric intersection (a stick).
\[\includegraphics[width=.2\textwidth]{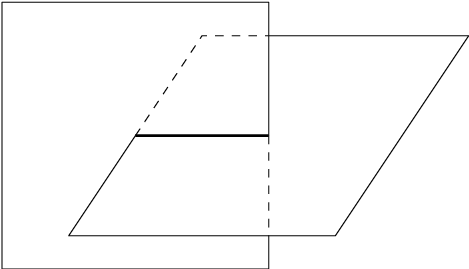}\hskip4em
\includegraphics[width=.2\textwidth]{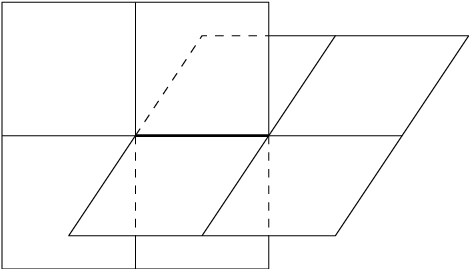}\]
\[\hbox{\small\sl Figure 10: (a) Rectangles in general position (b) Subdivided rectangles }\]

\medskip\noindent
Such rectangles can be decomposed into smaller rectangles (all with matching orientation) as in Figure 10(b). The product of a sum of four rectangles with another sum of four rectangles, thus splits by linearity into a sum of 16 intersections, four of which vanish due to being geometrically disjoint, leaving a sum of 12 non-zero terms, all of which are transverse pairs not in general position.
\[\includegraphics[width=.2\textwidth]{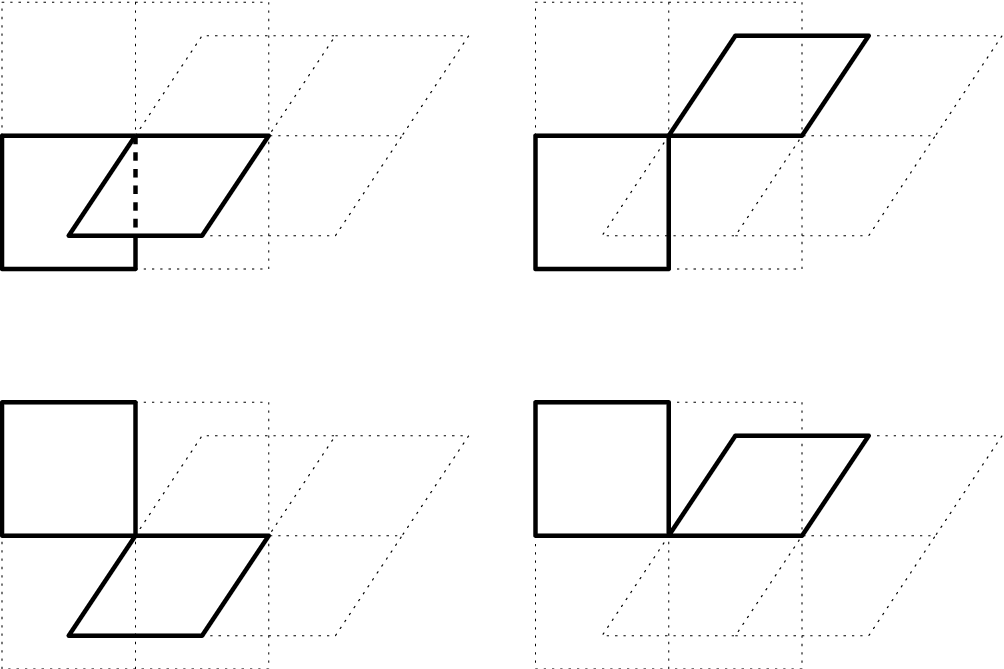}\hskip4em
\includegraphics[width=.2\textwidth]{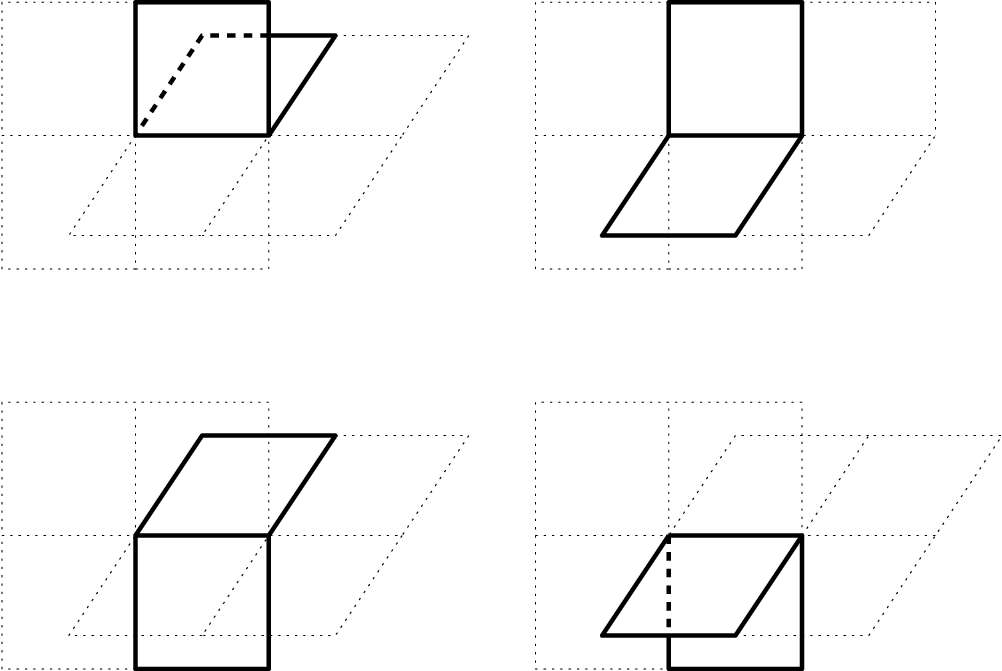}\hskip4em
\includegraphics[width=.2\textwidth]{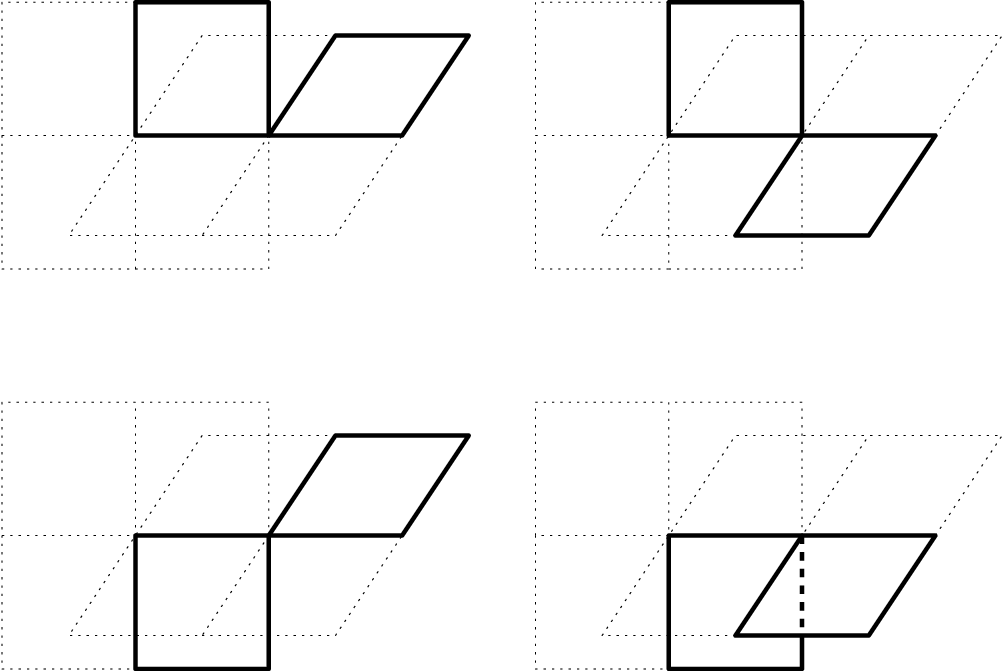}\]
\[\hbox{\small\sl Figure 11: 12 products comprising an intersection of rectangles in general position}\]

\medskip\noindent The left hand group of four products in Figure 11 all lead to an infinitesimal stick. All four figures can be obtained from each other via reflections in horizontal and vertical planes and thus by (D) have identical products. Similarly for the right hand set of four products (although they are located at a different point from the four products on the left). The middle group of products also for symmetry reasons consists of four identical terms, but each has support which is an $h$-stick.  Define the infinitesimal stick by
\[\includegraphics[width=.2\textwidth]{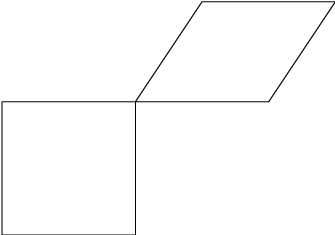}\hskip1ex\raise5ex\hbox{product$=\frac14$}
\includegraphics[width=.018\textwidth]{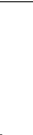}\eqno{(2)}\]
This can be done since this product is a multiple of an infinitesimal stick (by locality) and therefore the infinitesimal stick generator can be rescaled if necessary so that the coefficient here is precisely $\frac14$.

\medskip
Now we return to the product in Figure 10(a) which  by (F), being an intersection in general position, is the $h$-stick along the geometric intersection. We have expressed this product as a sum of 12 terms as in Figure 8 above. The sum of the four terms on the left is an infinitesimal stick at the left hand endpoint of the $h$-stick with weighting $4\cdot\frac14=1$ while the sum of the four terms on the right is an infinitesimal stick at the right hand endpoint of the $h$-stick (similarly weighted by 1). Therefore the sum of the remaining four terms in the middle must be the difference, the $h$-stick minus two infinitesimal sticks at the two endpoints. Since the four middle terms are equal (by (D)) thus each must be a quarter of this. Together this gives all non-zero products of $h$-squares.
\[\includegraphics[width=.4\textwidth]{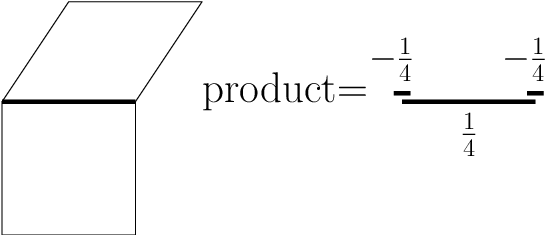}\eqno{(3)}\]

\medskip Concerning the products of squares and sticks $C_2\times{}FC_1\to{}C_0$, there are two types, namely $C_2\times{}C_1\to{}C_0$ which only involves the original $h$-complex $C_*$, and those products involving infinitesimal sticks. An intersection of a two-dimensional cuboid and a one-dimensional cuboid in general position looks like Figure 12(a). Property (F) requires that their product be the natural geometric intersection, namely a point (with weighting 1). 
\[\includegraphics[width=.1\textwidth]{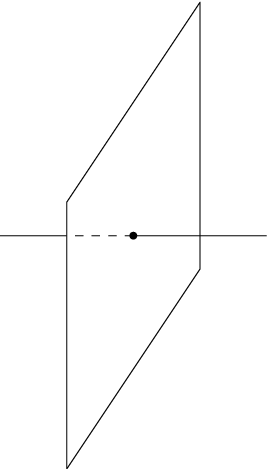}\hskip4em
\includegraphics[width=.1\textwidth]{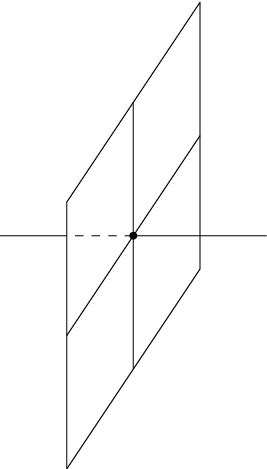}\]
\[\hbox{\small\sl Figure 12:(a) General position intersection of square and stick (b) Subdivided cells}\]

\medskip\noindent
Take the case of a $2h$-square and $2h$-stick with the configuration of Figure 12(a). Decomposing the square into four basis cells in the $h$-complex (all with matching orientation) and the stick into two $h$-sticks, as in Figure 12(b), then the product splits by linearity into a sum of eight products, each of an $h$-square with an orthogonal $h$-stick intersecting at a vertex. Each of these eight products can be obtained from each other by suitable symmetries of the lattice, and therefore by (D) they are all equal. We deduce that a non-zero product of an $h$-square and $h$-stick (which requires a point geometric intersection and the plane of the square to be orthogonal to the stick) is always, up to sign (dependent on the relative orientations of the square and stick) an eighth of the point 
\[\includegraphics[width=.1\textwidth]{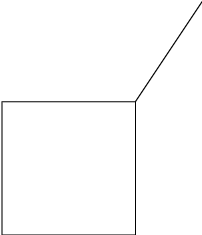}\hskip1ex\raise5ex\hbox{product$=\frac18\cdot$}\eqno{(4)}\]
Finally for the product of a $h$-square and an infinitesimal stick, which for a non-zero product must be directed orthogonally to the square and located at one of the vertices of the square. Since such an orthogonal stick arises as (four times) the intersection of a suitable pair of $h$-squares with only point intersection, we can use associativity to deduce its value
\[\includegraphics[width=.8\textwidth]{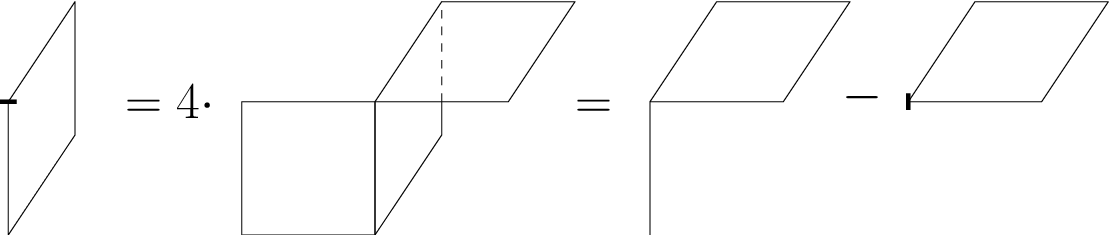}\]
By above the first term on the right hand side is an eighth of a point and thus the second term, which, being the multiple of a point must be by (D) the same as the left hand side, must be a sixteenth of a point,
\[\includegraphics[width=.06\textwidth]{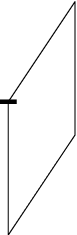}\hskip1em\raise7ex\hbox{product$=\frac{1}{16}\cdot$}\eqno{(5)}\]
Hence we have proved uniqueness, having identified in (2), (3), (4), (5) all non-trivial products of basis elements in $FC_*$. This concludes the proof, existence already having been demonstrated in \S4.
\end{proof}

\section{Addendum}
For arbitrary dimension $d$, the $d$-fold tensor product of the one-dimensional algebra with ideal sticks from \S3 gives an algebra $EC_*$ which is an extension of the cubical complex $C_*$ with ideal elements added 
  in all dimensions $1,\ldots,d$ and whose (transverse intersection) product satisfies properties (A)-(G) of the Comprehensive Theorem (p5, \S2). The part of $EC_*$ generated as an algebra by elements of the cubical complex with dimensions $\leq{}m$ defines a subalgebra $FC_*$ in which there are ideal elements only in dimensions $\leq2m-n$. The same argument as used in the proof of (H) of the Comprehensive Theorem shows that the product on $FC_*$ which is   graded commutative and  associative  also  satisfies the product rule everywhere, so long as $3m\leq2n$. 
 We use  the point of that proof  in   that this truncation limit   will mean products involving ideal elements will only have intersections in at most degree zero and so the product rule holds  because both sides are zero.  Otherwise there is trouble. So  $m=\left[\frac{2n}3\right]$ is the maximal truncation to obtain  commutative, associative and the product rule for the boundary operator on the  commutative and associative subalgebra generated by the truncated  $C_*$.

  In particular, {\bf in four dimensions}  to have the three properties commutative associative and product rule it is necessary to truncate to $C_{\leq2}$;  then there will be no ideal elements, $FC_*=C_{\leq2}$ in this case.

  {\bf In five dimensions}, $m=3$ and truncation to dimension at most three, $C_{\leq3}$ generates ideal elements in dimension one only (ideal sticks) and this  subalgebra has the three properties. 

  {\bf In six dimensions}, $m=4$ and $C_{\leq4}$ generates ideal elements in dimensions one and two and the subalgebra generated  has the three properties. Note there is an interesting 
 symmetric cubic three form generated by triple intersection followed by augmentation (summing coefficients in degree zero).


\begin{thebibliography}{100}
\bibitem{AnKwon} D.~An, A.~Kwon, R.~Lawrence, Probabilistic intersection and its application in fluid dynamics, {\sl draft report} (2023) 

\bibitem{LRS} R.~Lawrence, N.~Ranade, D.~Sullivan, Quantitative towers in finite difference calculus approximating the continuum, {\sl Quart.~J.~Math.} {\bf 72} (2021) 515--545

\bibitem{S} D.~ Sullivan, Lattice Hydrodynamics, Jean-Christoph\'e Yoccoz Memorial Volume, {\sl Ast\'erisque} {\bf 415} (2020) 215--222
	\end{thebibliography}
\end{document}